\documentclass[graybox]{svmult}

\usepackage{mathptmx}       % selects Times Roman as basic font
\usepackage{helvet}         % selects Helvetica as sans-serif font
\usepackage{courier}        % selects Courier as typewriter font
\usepackage{type1cm}        % activate if the above 3 fonts are
\usepackage{makeidx}         % allows index generation
\usepackage{graphicx}        % standard LaTeX graphics tool
\usepackage{multicol}        % used for the two-column index
\usepackage[bottom]{footmisc}% places footnotes at page bottom

\makeindex             % used for the subject index
\usepackage{amsmath}
\usepackage{amsfonts}
\newcommand{\Ex}{\mathbb{E}}
\newcommand{\Pb}{\mathbb{P}}
\newcommand{\Pn}{\mathbb{P}_n}
\newcommand{\Gn}{\mathbb{G}_n}
\newcommand{\Gb}{\mathbb{G}}
\newcommand{\R}{\mathbb{R}}
\def\ind{\perp\!\!\!\!\perp}

\newcommand{\eff}{{\scriptsize{\text{eff}}}}

\def\inprob{\stackrel{p}{\rightarrow}}
\def\indist{\rightsquigarrow}

\def\T{{ \mathrm{\scriptscriptstyle T} }}

\begin{document}

\title*{Semiparametric Theory and Empirical Processes in Causal Inference}
\author{Edward H. Kennedy}
\institute{Edward H. Kennedy \at University of Pennsylvania, Philadelphia, PA, USA, \email{kennedye@mail.med.upenn.edu}}
%\and Jason A. Roy \at University of Pennsylvania, Philadelphia, PA, USA, \email{jaroy@mail.med.upenn.edu}}

\maketitle

\vspace*{-.3in} 
\abstract{ In this paper we review important aspects of semiparametric theory and empirical processes that arise in causal inference problems. We begin with a brief introduction to the general problem of causal inference, and go on to discuss estimation and inference for causal effects under semiparametric models, which allow parts of the data-generating process to be unrestricted if they are not of particular interest (i.e., nuisance functions). These models are very useful in causal problems because the outcome process is often complex and difficult to model, and there may only be information available about the treatment process (at best). Semiparametric theory gives a framework for benchmarking efficiency and constructing estimators in such settings. In the second part of the paper we discuss empirical process theory, which provides powerful tools for understanding the asymptotic behavior of semiparametric estimators that depend on flexible nonparametric estimators of nuisance functions. These tools are crucial for incorporating machine learning and other modern methods into causal inference analyses. We conclude by examining related extensions and future directions for work in semiparametric causal inference. } %abstract word count = 178

\section{Introduction}
\label{sec:intro}

Causality and counterfactual questions lie at the heart of many if not most scientific endeavors. Counterfactual questions are about what \textit{would have happened} in some system had it undergone a particular change. For example: How would the distribution of patient outcomes differ had everyone versus no one received some medical treatment? Which rule for treatment assignment would maximize outcomes if it were implemented in the population?

In fact many scientific questions are causal even if they are not framed using explicitly causal language and notation. For example, standard regression analyses are often explained in implicitly causal terms, e.g., when regression coefficients are portrayed as representing the expected difference in outcome if all covariates were held constant, except for one covariate whose value was increased by one. In contrast, without causal assumptions, these coefficients can only represent the expected difference in outcome for two units who happen to have the same covariate values, except for one covariate whose values happen to differ by one; manipulation of the covariate cannot be allowed without invoking causal assumptions. 

In this chapter we give a review of semiparametric theory and empirical processes as they arise in causal inference problems. These include very powerful methodological tools that can be especially useful in causal settings. %In the rest of this introduction we describe the general layout of the chapter.

In Section \ref{sec:setup} we give an introduction to causal inference, following Robins \cite{robins1986new,robins2000marginal,vanderlaan2003unified}, van der Laan \cite{vanderlaan2003unified,vanderlaan2011targeted,vanderlaan2014targeted}, and others. In order to answer causal questions with observed data, we need causal assumptions. Sometimes these causal assumptions can hold by virtue of the study design (e.g., in randomized trials), while at other times the assumptions we need are untestable and need to be justified based on subject matter expertise (e.g., in standard observational studies). In either case, as we discuss in detail in Section \ref{subsec:semi_param}, it is important to have a clearly defined study question (with a corresponding causal parameter of interest). It is similarly important to be precise about the assumptions that are required to estimate the causal parameter of interest with observed data. This is the enterprise of identification, which we discuss briefly in Section \ref{subsec:semi_id}.

After a causal parameter of interest has been precisely defined and identified (i.e., expressed in terms of observed data), then  estimation and inference for that parameter is essentially a purely statistical problem. Classical maximum likelihood approaches can in theory be used to estimate such identified causal parameters, but typically require unrealistic parametric assumptions about the entire data-generating process. In contrast, semiparametric methods allow parts of the data-generating process to be completely unrestricted, e.g., if they are unknown or involve nuisance functions that are not of particular interest to the study question. Thus, if investigators have a good understanding of the treatment assignment process, for example, this information can be incorporated into a semiparametric analysis, and no assumptions might be needed about the outcome process. This is particularly useful in causal inference settings since the outcome process is often complex and difficult to model, while investigators may have some information about the treatment mechanism (e.g., by surveying doctors about how they prescribe some treatment). 

Alternatively, in many cases investigators may not have much information available about any part of the data-generating process. Then it will often be most reasonable to use a nonparametric model, which does not make any parametric assumptions at all about the data-generating process. A nonparametric model can be viewed as a special case of a semiparametric model, so the theory reviewed in this chapter covers these settings as well as those where treatment is assigned according to some known process.

In Section \ref{sec:semi} we review semiparametric theory, following foundational work by numerous authors, including Begun et al. \cite{begun1983information}, Bickel et al. \cite{bickel1993efficient}, Pfanzagl \cite{pfanzagl1982contributions}, van der Vaart \cite{vandervaart2000asymptotic,vandervaart2002semiparametric}, Robins \cite{robins1986new,robins2000marginal,vanderlaan2003unified}, van der Laan \cite{vanderlaan2003unified,vanderlaan2011targeted,vanderlaan2014targeted}, and many others \cite{tsiatis2006semiparametric,kosorok2007intro}. We start in Section \ref{subsec:semi_models} with a general introduction to semiparametric models, and discuss influence functions as representations of estimators in such models in Section \ref{subsec:semi_influence}. Then in Section \ref{subsec:semi_tangent} we introduce the notion of tangent spaces and a related space where influence functions reside, give an example illustrating basic semiparametric theory for estimation of the average treatment effect in Section \ref{subsec:eif}, and wrap up by discussing links to general missing data problems in Section \ref{subsec:full}.

Semiparametric theory gives us efficiency benchmarks in models where parts of the data-generating process are unrestricted, and tells us how to construct potentially efficient estimators. However, in order to understand the asymptotic behavior of such semiparametric estimators, particularly when flexible nonparametric methods are used to estimate nuisance functions, we need empirical process theory. This is the topic of Section \ref{sec:empirical}. The field of empirical processes is vast, so we only discuss parts that especially relate to estimation of nuisance functions. Our review follows important work by Andrews \cite{andrews1994empirical, andrews1994asymptotics}, Pollard \cite{pollard1984conv, pollard1990empirical}, van der Vaart \cite{vandervaart1996weak, vandervaart2000asymptotic, vandervaart2002semiparametric}, Wellner \cite{shorack1986empirical, vandervaart1996weak}, and others \cite{kosorok2007intro, vanderlaan2011targeted}. We start by giving the motivation for empirical process theory in semiparametric problems in Section \ref{subsec:motiv}, discuss Donsker classes and examples in Sections \ref{subsec:donsker} and \ref{subsec:donsker_ex}, and illustrate with an analysis of the doubly robust estimator of the average treatment effect in Section \ref{subsec:ate_ex}. 

We close the chapter in Section \ref{sec:future} by considering extensions and future directions for work in semiparametric causal inference.

\section{Setup}
\label{sec:setup}

In this section we briefly introduce the basic setup of a typical causal inference problem. We focus on two essential components of causal inference: first, formulating a clearly defined parameter of interest, and second, exploring how and whether this target parameter is identified with observed data. These issues are very important and provide a crucial foundation for semiparametric causal inference; however, we give only a brief treatment since the main goal of this chapter is to discuss semiparametric theory and empirical processes. Much of the discussion here is inspired by pioneering work by Robins \cite{robins1986new,robins2000marginal,vanderlaan2003unified}, van der Laan \cite{vanderlaan2003unified,vanderlaan2011targeted,vanderlaan2014targeted}, and colleagues.

\subsection{The Target Parameter}
\label{subsec:semi_param}

An important first step in any scientific pursuit is to have a clearly defined goal. In a statistical analysis, this includes giving a precise expression for a parameter of interest, which we will refer to as \textit{the target parameter}. 

The target parameter is the main feature of interest in the analysis, and ideally is decided upon based on collaborative discussion between scientific investigators and the statistician or analyst. In practice, however, the target parameter is sometimes defined only in vague terms, or is chosen based on convenience rather than scientific interest. In causal inference problems, the target parameter is typically formulated in terms of hypothetical interventions and corresponding {counterfactual} data, which represent the data that {would have} been observed under some intervention. In this chapter we mostly rely on the potential outcome framework, due to Neyman \cite{neyman1923sur} and Rubin \cite{rubin1974estimating,rubin1978role}, but note that alternative frameworks based on structural equation models and graphs \cite{pearl1995causal,pearl2009causality}, or decision theory \cite{dawid2000causal} can also be useful. 

For example, in some population of units (e.g., patients), let $Y \in \R$ denote a random variable representing an outcome of interest (e.g., blood pressure, or an indicator for whether a heart attack occurred), and let $A \in \{0,1\}$ denote a binary treatment (e.g., receipt of a statin), whose effect is in question. Then it may be of interest to estimate the average causal effect, i.e., how the expected outcome would have differed had everyone in the population taken treatment versus if no one in the population had taken treatment. This quantity can be represented notationally as follows. Let $Y^a$ denote the potential outcome that {would have} been observed (for a particular unit in the population) had that unit taken treatment level $A=a$. For a binary treatment, for example, this notation gives rise to two potential outcomes, $Y^1$ and $Y^0$, which are the outcomes that would have been observed for a particular unit under treatment ($A=1$) and control ($A=0$), respectively. Then the \textit{average causal effect} in the population can be defined as
\begin{equation}
\psi = \Ex(Y^1 - Y^0).
\label{eq:ace}
\end{equation}
Of course, different contrasts may instead be of interest under this hypothetical intervention; for example, if the outcome is binary then one may be more concerned with the risk ratio $\Ex(Y^1)/\Ex(Y^0)=\Pb(Y^1=1)/\Pb(Y^0=1)$, or with the odds ratio $\{\Pb(Y^1=1)/\Pb(Y^1=0)\}/\{\Pb(Y^0=1)/\Pb(Y^0=0)\}$. Alternatively, one may care more about how the effect of treatment changes with some other variable. Or some other entirely different intervention may be of interest; for example, one may want to learn what the mean outcome would have been if treatment had been assigned via some rule based on other variables \cite{murphy2003optimal,chakraborty2013statistical}, or how outcomes would have changed under treatment versus control if a mediating variable (a variable occurring subsequent to treatment, but prior to outcome) was fixed at some value \cite{tchetgen2012semiparametric, vanderweele2015explanation}. 

We will consider a number of different types of causal parameters and hypothetical interventions in subsequent sections, but a full taxonomy is beyond the scope of this chapter. The main point is that it is necessary to have a clear definition of the target parameter (i.e., the object one wants to learn about using data) when working in the semiparametric framework. In fact, regardless of framework or philosophical perspective, a clearly defined target parameter is necessary in order to meaningfully address estimation bias or variance relative to any meaningful standard.

\subsection{Identification}
\label{subsec:semi_id}

Once a target parameter is clearly defined based on some hypothetical intervention, the next step is to explore how and whether it can be \textit{identified} (i.e., expressed uniquely in terms of a distribution for observed data). This step translates the causal question of interest into a statistical problem defined in terms of observed data. 

For example, suppose that in a population of interest we actually get to observe potential outcomes under the received treatment for each unit, i.e., 
\begin{equation}
A=a \implies Y=Y^a . \tag{C1}
\label{eq:cons}
\end{equation}
Condition \eqref{eq:cons} is called ``consistency'' \cite{vanderweele2009concerning} and holds if potential outcomes are defined uniquely by a unit's own treatment and not others' (i.e., no interference), and also not by the way treatment is administered (i.e., no different versions of treatment). Also suppose that there exists some set of observed covariates $L$ that render treatment independent of potential outcomes when conditioned upon, i.e.,
\begin{equation}
A \ind Y^a \mid L  ,  \tag{C2}
\label{eq:ignor}
\end{equation}
where $\ind$ denotes statistical independence. Condition \eqref{eq:ignor} is often called ``no unmeasured confounding'', ``exchangeability'', or ``ignorability'', and holds if treatment is externally randomized, or if treatment decisions are made based only on covariates $L$. Finally suppose that, regardless of covariate value, each unit has a non-zero chance to receive treatment level $A=a$, i.e.,
\begin{equation}
p(A=a \mid L=l) \geq \delta > 0 \ \text{ whenever } p(L=l)>0 ,  \tag{C3}
\label{eq:pos}
\end{equation}
where $p(\cdot)$ denotes densities with respect to an appropriate dominating measure. Condition \eqref{eq:pos} is called ``positivity'' and means treatment is not assigned deterministically \cite{petersen2010diagnosing}. Then, if Conditions \eqref{eq:cons}--\eqref{eq:pos} hold for treatment value $a$, it follows that
\begin{equation}
p(Y^a=y \mid L=l) = p(Y=y \mid L=l, A=a) .
\end{equation}
Therefore we can express the conditional distribution of the potential outcome $Y^a$ given $L$ in terms of observed data; thus we can also identify the conditional distribution given any subset of $L$, including the null set, by simply marginalizing. In particular if Conditions \eqref{eq:cons}--\eqref{eq:pos} hold for $a=0,1$, then the average causal effect $\psi$ from \eqref{eq:ace} can be written as
\begin{equation}
\psi = \int_\mathcal{L} \Big\{ \Ex(Y \mid L=l, A=1) - \Ex(Y \mid L=l, A=0) \Big\} \ dP(L=l) .
\label{eq:ace2}
\end{equation}

The above identification result is an example of the g-computation formula, which was first proposed for general time-varying treatments by Robins \cite{robins1986new, robins2009estimation}. Numerous alternative identification schemes are also available, for example based on instrumental variables \cite{angrist1996ident,hernan2006instruments}. The literature on causal identification is extensive, and includes graphical criteria \cite{pearl1995causal, pearl2009causality}, bounds \cite{manski2003partial}, and many other topics. 

In this chapter we focus on settings where the target causal parameter (call it $\psi$) is identified, and thus can be written in terms of the distribution $P$ of the observed data. In the next section we illustrate ideas with the average causal effect $\psi$ defined in Equation \eqref{eq:ace}, and defined by Equation \eqref{eq:ace2} under Conditions \eqref{eq:cons}--\eqref{eq:pos}; although we focus on simple average effects, the general logic is similar for other parameters.

\section{Semiparametric Theory}
\label{sec:semi}
% Always give a unique label
% and use \ref{<label>} for cross-references
% and \cite{<label>} for bibliographic references
% use \sectionmark{}
% to alter or adjust the section heading in the running head

In this section we give a general review of semiparametric theory, using as a running example the common problem of estimating an average causal effect. Our review draws on foundational work in general semiparametric theory by Begun et al. \cite{begun1983information}, Bickel et al. \cite{bickel1993efficient}, Pfanzagl \cite{pfanzagl1982contributions, pfanzagl1990estimation}, and van der Vaart \cite{vandervaart2000asymptotic,vandervaart2002semiparametric}, among others \cite{newey1994asymptotic, kosorok2007intro}, as well as further developments for missing data and causal inference problems by Robins \cite{robins1986new, robins1994estimation, robins1995analysis,  robins2000marginal, vanderlaan2003unified}, van der Laan \cite{vanderlaan2003unified, vanderlaan2011targeted, vanderlaan2014targeted}, and colleagues \cite{hahn1998role, tsiatis2006semiparametric}.

\subsection{Semiparametric Models}
\label{subsec:semi_models}

Standard semiparametric theory generally considers the following setting. We observe an independent and identically distributed sample $(Z_1,...,Z_n)$ distributed according to some unknown probability distribution $P_0$ on the Borel $\sigma$-field $\mathcal{B}$ for some sample space $\mathcal{Z}$. The general goal is estimation and inference for some target parameter $\psi_0=\psi(P_0) \in \R^p$, where $\psi=\psi(P)$ can be viewed as a map from a probability distribution to the parameter space (assumed to be Euclidean here). In our running example where $\psi$ is the average causal effect defined in \eqref{eq:ace2} (after imposing identifying assumptions), the observed data consist of an independent and identically distributed sample of $Z=(L,A,Y)$ where $L$ denotes covariates, $A$ is a binary treatment, and $Y$ is the outcome of interest. Here we suppose the distribution $P_0$ has density given by
\begin{equation}
p(z) = p(y \mid l, a) p(a \mid l) p(l)
\label{eq:dens}
\end{equation}
with respect to some dominating measure. In general we write $p(X=t)$ for the density of $X$ at $t$, but when there is no ambiguity we let $p(x)=p(X=x)$.  

A \textit{statistical model} $\mathcal{P}$ is a set of possible probability distributions, which is assumed to contain the observed data distribution $P_0$. In a parametric model, $\mathcal{P}$ is assumed to be indexed by a finite-dimensional real-valued parameter $\theta \in \R^q$, e.g., we may have $\mathcal{P}=\{ P_\theta : \theta \in \R^q\}$ with $\psi \subseteq \theta$. For example, if $Z$ is a scalar random variable one might assume it is normally distributed with unknown mean and variance, $Z \sim N(\mu,\sigma^2)$, in which case the model is indexed by $\theta=(\mu,\sigma^2) \in \R \times \R^+$. \textit{Semiparametric models} are simply sets of probability distributions that cannot be indexed by only a Euclidean parameter, i.e., models that are indexed by an infinite-dimensional parameter. Semiparametric models can vary widely in the amount of structure they impose; for example, they can range from \textit{nonparametric models} for which $\mathcal{P}$ consists of all possible probability distributions, to simple regression models that characterize the regression function parametrically but leave the residual error distribution unspecified. 

In semiparametric causal inference settings it is common to impose some structure on the treatment mechanism (e.g., with a parametric model) leaving the outcome mechanism unspecified. This is because the outcome mechanism is often a complex natural process outside of human control, whereas the treatment mechanism is known in randomized trials, and can be well-understood in some observational settings (for example, when a medical treatment is assigned in a standardized way, which is communicated by physicians to researchers). In our running example, one may wish to do inference for the average causal effect $\psi$ under a parametric model for the treatment mechanism, leaving everything else unspecified, so that
\begin{equation}
p(z;\eta,\alpha) = p(y \mid l, a; \eta_y) p(a \mid l; \alpha) p(l;\eta_l) ,
\end{equation}
where $\alpha \in \R^q$ but $\eta=(\eta_y,\eta_l)$ represents an infinite-dimensional parameter that does not restrict the distribution of the outcome given covariates and treatment $p(y \mid l, a)$ or the marginal covariate distribution $p(l)$. 

Of course it is not always the case that there is substantive information available about the treatment mechanism; in many observational studies neither the exposure nor the outcome process is under human control, and both processes may be equally complex (e.g., in studies where the treatment or exposure is itself a disease or other medical condition). In such cases it is often more appropriate to consider inference for $\psi$ under a nonparametric model that makes no parametric assumptions about the distribution $P$. As we will see in Section \ref{subsec:ate_ex}, in order to obtain usual root-n rates of convergence in nonparametric models, we will still require some conditions on how well we can estimate the nuisance functions.

Another way semiparametric models arise in causal settings is through parametric assumptions about high-level treatment effects. For example, suppose we were not interested in the average causal effect $\Ex(Y^1-Y^0)$ but in how this effect varied with a subset of covariates $V \subset L$, i.e., the goal was to estimate $\gamma(v)=\Ex(Y^1 - Y^0 \mid V=v)$. Letting $W=L \setminus V$ so that $L=(V,W)$, it is straightforward to show that this conditional effect is also identified under Conditions (C1)--(C3) as in \eqref{eq:ace2}, except replacing $dP(l)$ with $dP(w \mid v)$. If $V$ includes a continuous variable or has many strata, it may be desirable to make parametric assumptions to reduce the dimension of $\gamma(v)$ (or in rare cases, there may be substantive knowledge about the parametric form of the effect modification), and thus one may want to assume $\gamma(v)=\gamma(v;\psi)$ for $\psi \in \R^p$. Such assumptions are not always easily encoded directly in the distribution $p(z)$, but can still be employed in conjunction with parametric assumptions about the treatment mechanism, for example, or in otherwise nonparametric models. An alternative approach is to use nonparametric \textit{working models} \cite{neugebauer2007nonparametric}, where instead of assuming $\gamma(v)=\gamma(v;\psi)$ we define our target parameter as a projection of $\gamma(v)$ onto the model $\gamma(v;\psi)$ (using, for example, a weighted least squares projection). %Examples of both approaches will be discussed in detail in Section \ref{sec:example}. 

\subsection{Influence Functions}
\label{subsec:semi_influence}

In the previous subsection we discussed the concept of a semiparametric model (in which part of the distribution $P$ is allowed to have unrestricted or infinite-dimensional components) and gave some examples. Now we begin to discuss estimation and inference in such models. This requires the concept of the \textit{influence function}, which is a foundational object of statistical theory that allows us to characterize a wide range of estimators and their efficiency.

Let $\Pn=n^{-1} \sum_i \delta_{Z_i}$ denote the empirical distribution of the data, where $\delta_z$ is the Dirac measure that simply indicates whether $Z=z$. This means for example that empirical averages can be written as $n^{-1} \sum_i f(Z_i) = \int f(z) \ \!d\Pn =\Pn\{ f(Z)\}$. An estimator $\hat\psi=\hat\psi(\Pn)$ is asymptotically linear with influence function $\varphi$ if the estimator can be approximated by an empirical average in the sense that
\begin{equation}
\hat\psi - \psi_0 = \Pn\{ \varphi(Z) \} + o_p(1/\sqrt{n}) ,
\label{eq:inf_fn}
\end{equation}
where $\varphi$ has mean zero and finite variance (i.e., $\Ex\{\varphi(Z)\}=0$ and $\Ex\{\varphi(Z)^{\otimes 2}\}< \infty$). Here $o_p(1/\sqrt{n})$ employs the usual stochastic order notation so that $X_n=o_p(1/r_n)$ means $r_n X_n \inprob 0$ where $\inprob$ denotes convergence in probability. 

Importantly, by the classical central limit theorem, an estimator $\hat\psi$ with influence function $\varphi$ is asymptotically normal with
\begin{equation}
\sqrt{n}(\hat\psi - \psi_0) \indist N\Big(0, \ \Ex\{\varphi(Z)^{\otimes 2}\} \Big) ,
\end{equation}
where $\indist$ denotes convergence in distribution. Thus if we know the influence function for an estimator, we know its asymptotic distribution, and we can easily construct confidence intervals and hypothesis tests, for example. Also, the efficient influence function for an asymptotically linear estimator is almost surely unique (i.e., unique up to measure zero sets) \cite{tsiatis2006semiparametric}, so in a sense the influence function contains all information about an estimator's asymptotic behavior (up to $o_p(1/\sqrt{n})$ error). 

Consider our running example where $\psi$ is the average causal effect defined in Equations \eqref{eq:ace} and \eqref{eq:ace2}. Suppose we are in a randomized trial setting where the propensity score $\pi(l) = p(A=1 \mid L=l)$ is known. A simple inverse-probability weighted estimator is given by
\begin{equation}
\hat\psi_{ipw} = \Pn\left\{ \frac{AY}{\pi(L)} - \frac{(1-A)Y}{1-\pi(L)} \right\} .
\label{eq:ipw}
\end{equation}
(Note that $\Ex(\hat\psi_{ipw})=\psi_0$ by iterated expectation.) The influence function for the estimator $\hat\psi_{ipw}$ is clearly given by 
\begin{equation}
\varphi_{ipw}(Z) = \frac{AY}{\pi(L)} - \frac{(1-A)Y}{1-\pi(L)} - \psi_0
\end{equation}
since $\hat\psi_{ipw} - \psi_0 = \Pn\{\varphi_{ipw}(Z)\}$ exactly, without any $o_p(1/\sqrt{n})$ approximation error. 

Now suppose we are in an observational study setting where the propensity score $\pi(l)$ needs to be estimated, and suppose we do so with a correctly specified parametric model $\pi(l;\alpha)$, with $\alpha \in \R^q$, so that the estimator $\hat\alpha$ solves some estimating equation $\Pn\{S(Z;\hat\alpha)\}=0$. Then the inverse-probability-weighted estimator $\hat\psi_{ipw}^*$ is given by \eqref{eq:ipw} above, except with the estimated propensity score $\pi(L;\hat\alpha)$ replacing the true propensity score $\pi(L)$. We can find the corresponding influence function by standard estimating equation techniques \cite{boos2002calculus}. Specifically, we have that $\hat\theta=(\hat\psi_{ipw}^*,\hat\alpha^\T)^\T$ solves $\Pn\{ m(Z;\hat\theta)\}=0$ where $m(z;\theta)=\{\varphi_{ipw}(Z;\psi,\alpha),S(Z;\alpha)^\T\}^\T$ are the stacked estimating equations for $\psi$ and $\alpha$, with the influence function for known propensity score given by $\varphi_{ipw}(Z;\psi,\alpha)={AY}/{\pi(L;\alpha)} - {(1-A)Y}/\{1-\pi(L;\alpha)\} - \psi$. Then under standard regularity conditions \cite{newey1994large, vandervaart2000asymptotic, tsiatis2006semiparametric} we have
\begin{equation}
\hat\theta - \theta_0 = \Pn\left[ \Ex\left\{ \frac{\partial m(Z;\theta_0)}{\partial \theta} \right\}^{-1} m(Z;\theta_0) \right] + o_p(1/\sqrt{n}) ,
\end{equation}
which after evaluating and rearranging implies that the influence function for $\hat\psi_{ipw}^*$ when the propensity score $\pi(l;\alpha)$ is estimated is
\begin{equation*}
\varphi_{ipw}^*(Z) = \varphi_{ipw}(Z;\psi_0,\alpha_0) - \Ex\left\{ \frac{\partial \varphi_{ipw}(Z;\psi_0,\alpha_0)}{\partial \alpha^\T} \right\} \Ex\left\{ \frac{\partial S(Z;\alpha_0)}{\partial \alpha} \right\}^{-1} S(Z;\alpha_0) .
\end{equation*} 

Surprisingly, even if the propensity score is known, it can be shown \cite{tsiatis2006semiparametric} that the inverse-probability-weighted estimator $\hat\psi_{ipw}^*$ based on an estimated propensity score is at least as efficient as the inverse-probability-weighted estimator $\hat\psi_{ipw}$ that uses the known propensity score. In other words, the variance of the influence function $\varphi_{ipw}^*(Z)$ is less than or equal to the variance of the influence function $\varphi_{ipw}(Z)$ for known propensity score. Thus the propensity score should be estimated from the data (according to a correct model, of course) even when it is known; discarding information can actually yield better efficiency.

So far we have seen that, given an estimator $\hat\psi$, we can learn about its asymptotic behavior by considering its influence function $\varphi(Z)$. But we can also use influence functions to find or construct estimators. Suppose we are given a candidate influence function $\varphi(Z;\psi,\eta)$ that depends on the target parameter $\psi$ as well as a nuisance parameter $\eta$ as in the previous examples. Then we can construct an estimator by solving the estimating equation $\Pn\{ \varphi(Z;\psi,\hat\eta) \} = 0$ in $\psi$, where $\hat\eta$ is some estimate of the nuisance parameter. Under standard regularity conditions, along with some additional conditions on the nuisance estimation, the corresponding estimator will itself be asymptotically linear with an influence function related to $\varphi(Z;\psi_0,\eta_0)$ depending on the form of the function $\varphi$ and how the nuisance parameter $\eta$ is estimated (as in the previous example). Other approaches for constructing estimators based on a particular influence function are also possible \cite{vanderlaan2006targeted, vanderlaan2011targeted}.

There is a deep connection between (asymptotically linear) estimators for a given model and the influence functions under that model. In some sense, if we know one then we know the other. Thus if we can find all the influence functions for a given model, we can characterize all asymptotically linear estimators for that model.

\subsection{Tangent Spaces}
\label{subsec:semi_tangent}

In this subsection we discuss the fundamental problem of how to find influence functions for a given semiparametric model, by characterizing the space in which influence functions reside. As noted previously, once we have solved this problem we can characterize valid estimators under our model. In particular, we can use influence functions to construct estimators and explore their efficiency.

To ease notation, consider the case where the target parameter is a scalar, so that $\psi \in \R$. As discussed in the previous subsection, influence functions $\varphi$ are functions of the observed data $Z$ with mean zero and finite variance. These influence functions reside in the Hilbert space $L_2(P)$ of measurable functions $g: \mathcal{Z} \rightarrow \R$ with $P g^2 = \int g^2 \ dP = \Ex\{g(Z)^2\} < \infty$, equipped with covariance inner product $\langle g_1,g_2 \rangle = P(g_1 g_2)$. The space of influence functions will be a subspace of this Hilbert space. A Hilbert space is a complete inner product space, and can be viewed as a generalization of usual Euclidean space; it provides a notion of distance and direction for spaces whose elements are potentially infinite-dimensional functions. 

A fundamentally important subspace of $L_2(P)$ in semiparametric problems is the \textit{tangent space}. First we will discuss the tangent spaces for parametric models. For parametric models indexed by real-valued parameter $\theta \in \R^{q+1}$, the tangent space $\mathcal{T}$ is defined as the linear subspace of $L_2(P)$ spanned by the score vector, i.e.,
\begin{equation}
\mathcal{T} = \{ b^\T S_\theta(Z;\theta_0) : b \in \R^{q+1} \}, 
\end{equation}
where $S_\theta(Z;\theta_0) = \partial \log p(z;\theta) / \partial \theta |_{\theta=\theta_0}$. If we can decompose $\theta=(\psi,\eta)$ then we can equivalently write $\mathcal{T} = \mathcal{T}_\psi \oplus \mathcal{T}_\eta$ for
\begin{equation}
\mathcal{T}_\psi = \{ b_1 S_\psi(Z;\theta_0) : b_1 \in \R \} \ , \ \mathcal{T}_\eta = \{ b_2^\T S_\eta(Z;\theta_0) : b_2 \in \R^q \}, 
\end{equation}
where $S_\psi(Z;\theta_0) = \partial \log p(z;\theta) / \partial \psi |_{\theta=\theta_0}$ is the score function for the target parameter, and similarly $S_\eta(Z;\theta_0) = \partial \log p(z;\theta) / \partial \eta |_{\theta=\theta_0}$ is the score for the nuisance parameter ($A \oplus B$ denotes the direct sum $A \oplus B=\{a + b : a \in A, b \in B\}$). In the above formulation, the space $\mathcal{T}_\eta$ is called the \textit{nuisance tangent space}. Influence functions for $\psi$ reside in the \textit{orthogonal complement of the nuisance tangent space}, denoted by $\mathcal{T}_\eta^\perp = \{ g \in L_2(P) : P(g h) = 0 \ \text{ for any } \ h \in \mathcal{T}_\eta \}$. In such parametric settings, this orthogonal space $\mathcal{T}_\eta^\perp$ can be written as 
\begin{align}
\mathcal{T}_\eta^\perp &= \{ g \in L_2(P) : g= h - \Pi(h \mid \mathcal{T}_\eta) , \ h \in L_2(P) \}  \\
&= \{ g \in L_2(P) : g= h - P(h S_\eta^\T) P(S_\eta S_\eta^\T)^{-1} S_\eta ,\  h \in L_2(P) \} , \nonumber
\end{align}
where $\Pi(g \mid S)$ denotes projections of $g$ on the space $S$, i.e., $P[h \{g-\Pi(g \mid S)\}]=0$ for all $h \in S$. The subspace of influence functions is the set of elements $\varphi \in \mathcal{T}_\eta^\perp$ that satisfy $P(\varphi S_\psi)=1$. The \textit{efficient influence function} is the influence function with the smallest covariance $P(\varphi^2)$, and is given by $\varphi_\eff=P(S_\eff^2)^{-1} S_\eff$, where $S_\eff$ is the \textit{efficient score}, given by $S_\eff=S_\psi-\Pi(S_\psi \mid \mathcal{T}_\eta)$.

Thus if we can characterize the nuisance tangent space and its orthogonal complement, then we can characterize influence functions. In fact, one can show that all regular asymptotically linear estimators have influence functions $\varphi$ that reside in $\mathcal{T}_\eta^\perp$ with $P(\varphi S_\psi)=1$, and conversely any element in this space corresponds to the influence function for some regular asymptotically linear estimator \cite{tsiatis2006semiparametric}. Thus characterizing the nuisance tangent space allows us to also characterize all regular asymptotically linear estimators. (Recall that a regular estimator is one whose limiting distribution is insensitive to local changes to the data generating process, as defined for example in \cite{vandervaart2000asymptotic, tsiatis2006semiparametric} and elsewhere.)

We have seen that in parametric models the tangent space is defined as the span of the score vector $S_\theta$. However, in semiparametric models, the nuisance parameter is infinite-dimensional and cannot be indexed by a real-valued parameter, so we cannot define scores in the usual way, since this requires differentiation with respect to the nuisance parameter. How can we extend the concept of the tangent space to semiparametric settings?

Constructing tangent spaces in semiparametric models requires a technical device called a \textit{parametric submodel}. A parametric submodel $\mathcal{P}_\epsilon$ indexed by real-valued parameter $\epsilon$ is a set of distributions contained in the larger model $\mathcal{P}$, which also contains the truth (i.e., $P_0 \in \mathcal{P}_\epsilon$); typically we have $\mathcal{P}_\epsilon = \{ P_\epsilon : \epsilon \in \R\}$ with $P_\epsilon|_{\epsilon=0} = P_0$. Thus a parametric submodel needs to respect the semiparametric model $\mathcal{P}$ and also needs to equal the true distribution at $\epsilon=0$. A typical example of a parametric submodel is given by
\begin{equation}
p_\epsilon(z) = p_0(z) \{ 1 + \epsilon g(z)\} ,
\label{eq:submodel}
\end{equation}
where $\Ex\{g(Z)\}=0$ and we have $\sup_z |g(z)|  < M$ and $|\epsilon| < 1/M$ so that $p_\epsilon(z) \geq 0$. We will often index the parametric submodel by the function $g$, and so let $P_\epsilon=P_{\epsilon,g}$. Note again that parametric submodels like the one above are a technical device for constructing tangent spaces and analyzing semiparametric models, rather than a usual model whose parameters we want to estimate from data (since $P_\epsilon$ depends on the true distribution $P_0$, it cannot be used as a model in the usual sense) \cite{tsiatis2006semiparametric}.

One intuition behind parametric submodels can be expressed in terms of efficiency bounds as follows \cite{vandervaart2000asymptotic}. First note that it is an easier problem to estimate $\psi$ under the parametric submodel $\mathcal{P}_\epsilon \in \mathcal{P}$ than it is to estimate $\psi$ under the entire (larger) semiparametric model $\mathcal{P}$. Therefore the efficiency bound under the larger model $\mathcal{P}$ must be larger than the efficiency bound under any parametric submodel. In fact we can define the efficiency bound for semiparametric models as the supremum of all such parametric submodel efficiency bounds.

Now that we have defined parametric submodels, how can they be used to construct tangent spaces? Just as the tangent space is defined as the linear span of the score vector in parametric models, in semiparametric models the tangent space $\mathcal{T}$ is defined as the (closure of the) linear span of scores of the parametric submodels. In other words, we first define scores on the parametric submodels $P_\epsilon$ with $S_\epsilon(z) = \partial \log p_\epsilon(z)/\partial \epsilon|_{\epsilon=0}$, and then construct parametric submodel tangent spaces as described earlier for standard parametric models, i.e., $\mathcal{T}_\epsilon = \{ b^\T S_\epsilon(Z) : b \in \R \}$. Note that for parametric submodels like the one defined in \eqref{eq:submodel} we have
\begin{equation}
S_\epsilon(z) =  g(z)/\{1+\epsilon g(z)\}|_{\epsilon=0} = g(z) ,
\end{equation}
so that the functions $g$ indexing the parametric submodels are set up to equal the parametric submodel scores. The closure $\mathcal{T}$ of the parametric submodel tangent spaces $\mathcal{T}_\epsilon$ is the minimal closed set that contains them; roughly speaking, $\mathcal{T}$ is the union of all the spaces $\mathcal{T}_\epsilon$ along with their limit points. Similarly, the nuisance tangent space $\mathcal{T}_\eta$ for a semiparametric model is the set of scores in $\mathcal{T}$ that do not vary the target parameter $\psi$, i.e.,
\begin{equation}
\mathcal{T}_\eta = \{ g \in \mathcal{T} : \partial \psi(P_{\epsilon,g})/\partial \epsilon|_{\epsilon=0} = 0 \}.
\end{equation}
Importantly, in nonparametric models the tangent space is the whole Hilbert space of mean zero functions. For more restrictive semiparametric models the tangent space will be a proper subspace. 

Now that we are equipped with definitions of tangent spaces and nuisance tangent spaces in semiparametric models, we can define influence functions, efficient influence functions, and efficient scores in much the same way we did before with parametric models. 

Specifically, the subspace of influence functions is the set of elements $\varphi \in \mathcal{T}_\eta^\perp$ that satisfy $P(\varphi S_\psi)=1$. The {efficient influence function} is the influence function with the smallest covariance $P(\varphi_\eff^2) \leq P(\varphi^2)$ for all $\varphi$; it is given by $\varphi_\eff=P(S_\eff^2)^{-1} S_\eff$, where $S_\eff$ is the {efficient score} defined as the projection of the score onto the tangent space, i.e., $S_\eff=\Pi(S_\psi \mid \mathcal{T}_\eta^\perp) = S_\psi-\Pi(S_\psi \mid \mathcal{T}_\eta)$ as before. The efficient influence function can also be defined as the projection of any influence function $\varphi$ onto the tangent space, $\varphi_\eff = \Pi(\varphi \mid \mathcal{T})$ for any influence function $\varphi$, as well as the pathwise derivative of the target parameter in the sense that $P(\varphi S_\epsilon)=\partial \psi(P_{\epsilon})/\partial \epsilon|_{\epsilon=0}$.

\subsection{Efficient Influence Function for Average Treatment Effect}
\label{subsec:eif}

As an illustration, return to our example involving the average treatment effect $\psi=\Ex(Y^1-Y^0)=\Ex\{\mu(L,1)-\mu(L,0)\}$, where we let $\mu(l,a)=\Ex(Y \mid L=l, A=a)$ denote the outcome regression function. Also let $\pi(l)=P(A=1 \mid L=l)$ denote the propensity score as before. In this subsection, we will show using the results from previous subsections that, under a nonparametric model where the distribution $P$ is unrestricted, the efficient influence function for $\psi$ is given by $\varphi(Z;\psi,\eta) = m_1(Z;\eta) - m_0(Z;\eta) - \psi$, where
\begin{equation}
%\varphi(Z;\psi,\eta) = \frac{A\{Y-\mu(L,1)\}}{\pi(L)} - \frac{(1-A)\{Y-\mu(L,1)\}}{1-\pi(L)} + \mu(L,1)-\mu(L,0)-\psi
m_a(Z;\eta) = m_a(Z;\pi,\mu) = \frac{I(A=a) \{Y - \mu(L,a)\}}{a\pi(L) + (1-a)\{1-\pi(L)\}} + \mu(L,a)
\label{eq:eif}
\end{equation}
with $\eta=(\pi,\mu)$ the nuisance function for this problem.

We will show this result by checking that the proposed efficient influence function $\varphi$ is a pathwise derivative in the sense that $\partial \psi(P_{\epsilon})/\partial \epsilon|_{\epsilon=0}=P(\varphi S_\epsilon)$. 

Here we let $p_\epsilon(z)=p(z;\epsilon)$ denote a parametric submodel with parameter $\epsilon \in \R$. For notational simplicity let $f_\epsilon'(t;0)=\{\partial f(z;\epsilon)/\partial \epsilon\}|_{\epsilon=0}$ for any function $f$ of $\epsilon$ and $z$, and also let $\ell(v \mid w;\epsilon)=\log p(v \mid w;\epsilon)$ for any partition $(V,W) \subseteq Z$, so that for example scores on the parametric submodels are denoted by $S_\epsilon(z)=\ell'_\epsilon(z;0)$. Then by definition from \eqref{eq:ace2} we have
\begin{equation}
\ell'_\epsilon(z;\epsilon) = \ell'_\epsilon(y \mid l, a;\epsilon) + \ell'_\epsilon(a \mid l;\epsilon) + \ell'_\epsilon(l;\epsilon) .
\end{equation}

First consider the term $\partial \psi(P_{\epsilon})/\partial \epsilon|_{\epsilon=0}=\psi'_\epsilon(0)$. By definition we have  $\psi = \int \int \{ y \ dP(y \mid l, a=1) - y \ dP(y \mid l, a=0) \} \ dP(l)$, 
%\begin{equation}
%\psi = \int \int \{ y \ dP(y \mid l, a=1) - y \ dP(y \mid l, a=0) \} \ dP(l)
%\end{equation}
so that
\begin{align}
\label{eq:psi_deriv}
\psi'_\epsilon(\epsilon) &= \int \int \{ y \ell'_\epsilon(y \mid l, a=1;\epsilon) \ dP(y \mid l, a=1; \epsilon)  \\
& \hspace{.6in } - y \ell'_\epsilon(y \mid l, a=0;\epsilon) \ dP(y \mid l, a=0; \epsilon) \} \ dP(l;\epsilon) \nonumber \\
&+ \int \int \{ y \ dP(y \mid l, a=1; \epsilon) - y \ dP(y \mid l, a=0; \epsilon) \} \ell'_\epsilon(l;\epsilon) \ dP(l;\epsilon) , \nonumber
\end{align}
where we used the fact that $dP'_\epsilon(v \mid w;\epsilon)=\ell'_\epsilon(v \mid w;\epsilon) dP(v \mid w;\epsilon)$. This follows since $\partial \log f(\epsilon)/\partial \epsilon = \{\partial f(\epsilon)/\partial \epsilon\}/f(\epsilon)$ for general functions $f$ by definition of the logarithmic derivative. Recall that when we evaluate the above at $\epsilon=0$, we have $dP(y \mid l, a;0)=dP(y \mid l,a)$ and $dP(l;0)=dP(l)$. 
 
Now consider the term $P(\varphi S_\epsilon)=\Ex\{ \varphi(Z;\psi,\eta) \ell'_\epsilon(Z;0)\}$, which equals
\begin{align}
\Ex\Big[ & \{ m_1(Z;\eta) - m_0(Z;\eta) - \psi \}\{\ell'_\epsilon(Y \mid L, A;0) + \ell'_\epsilon(A \mid L;0) + \ell'_\epsilon(L;0)\} \Big] \nonumber \\
&= \Ex\left[\left\{ \frac{A}{\pi(L)} - \frac{1-A}{1-\pi(L)} \right\} Y \ell'_\epsilon(Y \mid L, A;0) + \{ \mu(L,1)-\mu(L,0) \} \ell'_\epsilon(L;0) \right] \nonumber \\  
&= \Ex \Big[ \Ex\{Y \ell'_\epsilon(Y \mid L, A=1;0) \mid L,A=1 \} - \Ex\{Y \ell'_\epsilon(Y \mid L, A=0;0) \mid L,A=0 \} \nonumber \\
& \hspace{.5in} +  \{ \mu(L,1)-\mu(L,0) \} \ell'_\epsilon(L;0) \Big] \nonumber  \\
&= \int \int \{ y \ell'_\epsilon(y \mid l, a=1; 0) \ dP(y \mid l, a=1) \label{eq:psi_cov} \\
& \hspace{.7in } - y \ell'_\epsilon(y \mid l, a=0;0) \ dP(y \mid l, a=0) \} \ dP(l) \nonumber \\
& \hspace{.5in } + \int \int \{ y \ dP(y \mid l, a=1) - y \ dP(y \mid l, a=0) \} \ell'_\epsilon(l;0) \ dP(l) .  \nonumber
\end{align}
The first equality follows from iterated expectation and the fact that, by usual properties of score functions, $\Ex\{\ell'_\epsilon(V \mid W;0) \mid W\}=0$. The second equality follows from iterated expectation, and the third follows by definition. 

Since the last expression for the covariance $P(\varphi S_\epsilon)$ in Equation \eqref{eq:psi_cov} equals the expression for $\psi'_\epsilon(\epsilon)$ from Equation \eqref{eq:psi_deriv} when evaluated at $\epsilon=0$, we have shown that $\varphi$ is in fact the efficient influence function.

\subsection{Full vs.\ Observed Data Influence Functions}
\label{subsec:full}

So far we have introduced the notion of a tangent space and discussed how influence functions $\varphi$ for regular asymptotically linear estimators can be viewed as elements of a subspace of the Hilbert space $L_2(P)$, namely the orthogonal complement of the nuisance tangent space, i.e., $\varphi \in \mathcal{T}_\eta^\perp$. We also illustrated how to check that a proposed influence function is the efficient influence function. But how does one find the space $\mathcal{T}_\eta^\perp$ in a given problem? In many cases this is a bit of an art: one conjectures the form of $\mathcal{T}_\eta^\perp$ and then checks that the conjectured space satisfies the required properties. For nonparametric models, one can sometimes deduce the form of the efficient influence function from the nonparametric maximum likelihood estimator, assuming discrete data \cite{vanderlaan2011targeted}. However, in some settings it can be useful to characterize influence functions with hypothetical `full data' (i.e., had we observed all counterfactuals), and then map these to observed data influence functions \cite{vanderlaan2003unified}. 

To characterize full-data influence functions in causal inference problems we need to start by presenting causal inference as a missing data problem \cite{vanderlaan2003unified, tsiatis2006semiparametric}. Thus far we have supposed that we observe an independent and identically distributed sample of observations $Z \sim P$. In general missing data problems, we conceive of hypothetical full data $\tilde{Z}$, of which the observed data $Z$ is a coarsened version. The problem is that we want to learn about the distribution $\tilde{P}$ of the full data $\tilde{Z}$, but we only get to observe the coarsened version $Z$ of the full data $\tilde{Z}$. In general coarsened data problems, $Z=\Phi(\tilde{Z},C)$ is a known many-to-one function $\Phi(\cdot)$ of both $\tilde{Z}$ and a coarsening variable $C$ that indicates what portion of $\tilde{Z}$ is observed. In causal inference settings, the coarsening variable generally equals the treatment process so that $C=A$, and
\begin{equation}
\tilde{Z} = \{Z^a : a \in \mathcal{A}\}.
\end{equation}
Thus the full data $\tilde{Z}$ are the potential outcomes under different levels $a \in \mathcal{A}$ of a general treatment process $A$ (here $A$ could be multivariate, e.g., a treatment sequence over multiple timepoints). For a given unit we only get to observe $Z=\Phi(\tilde{Z},A)=Z^A$, i.e., the potential outcome under the observed treatment process. For instance, in our running example where $Z=(L,A,Y)$ with binary treatment so that $\mathcal{A}=\{0,1\}$, the full data for a given unit could be represented as
\begin{equation}
\tilde{Z} = \{ (L^a,Y^a)  : a \in \{0,1\} \} = (L, Y^0, Y^1) .
\end{equation}
Note that the last equality follows since $L^a = L$ if we make the usual assumption that events in the past cannot be affected by the future. In some cases we might also want to include the observed treatment process in the full data, so that in the above example we would have $\tilde{Z}=(L,A,Y^0, Y^1)$. In a longitudinal setting where covariates and a binary treatment are updated at timepoints $t=1,...,K$ and an outcome is measured at the end of follow-up, we could have
\begin{equation}
\tilde{Z} = \{ (L_1, L_2^{a_1}, L_3^{a_1,a_2}, ..., L_t^{\overline{a}_{t-1}}, ..., L_K^{\overline{a}_{K-1}}, Y^{\overline{a}_K}) : \overline{a}_K \in \{0,1\}^K \},
\end{equation}
where $\overline{a}_t = (a_1,...,a_t)$ denotes the past history of a variable through time $t$. The observed data in this case would be $Z=(L_1,A_1,...,L_t,A_t,...,L_K,A_K,Y)$ for a given unit. Not every causal inference problem fits in the above framework, but when the framework applies it can often be very useful.

Now that we have defined the full data $\tilde{Z}$ and given some examples, we can also define corresponding tangent spaces, influence functions, and parametric submodels, using semiparametric models $\widetilde{\mathcal{P}}$ for the full data just as we did for the observed data previously. The advantage is that it is often more straightforward to derive tangent spaces and influence functions for full data problems (or else results may already be known for common models), and then translate them to observed data, rather than working with observed data directly and using the results from previous subsections. Of course, in order to translate full data influence functions to observed data influence functions, we need identifying assumptions. 

Under a coarsening at random assumption \cite{gill1997car}, results for mapping full data to observed data tangent spaces are given for example in \cite{vanderlaan2003unified} and \cite{tsiatis2006semiparametric}. In general, coarsening at random means $P(Z=z \mid \tilde{Z}=\tilde{z}_1) = P(Z=z \mid \tilde{Z}=\tilde{z}_2) \ \text{ whenever } \ z=\Phi(\tilde{z}_1,a)=\Phi(\tilde{z}_2,a)$ for some $a \in \mathcal{A}$. In many problems \cite{robins1999sensitivity}, this can be equivalently expressed by saying that $P(A=a \mid \tilde{Z}=\tilde{z}_1) = P(A=a \mid \tilde{Z}=\tilde{z}_2)$ only depends on $z$ whenever $z=\Phi(\tilde{z}_1,a)=\Phi(\tilde{z}_2,a)$. Under some conditions, coarsening at random also reduces to a  randomization assumption, which says treatment is independent of potential outcomes given the observed past, e.g., $A \ind Y^a \mid L$ in our running example, or $A_t \ind Y^{\overline{a}_K} \mid \overline{L}_t, \overline{A}_{t-1}$ in the above longitudinal example. More details on these issues are given in \cite{robins1999sensitivity, vanderlaan2003unified}. Again we point out that this framework does not always apply: sometimes coarsening at random is not equivalent to treatment randomization, or is not the identifying assumption we wish to utilize.

Here we will be content giving a simple example of how to map a full data influence function to the observed data, rather than discussing details in full generality; see \cite{vanderlaan2003unified} and \cite{tsiatis2006semiparametric} for more general results. Assume coarsening at random holds, and that the treatment assigment process is known. Further suppose the observed data is $Z=(L,A,Y)$ with $A \in \{0,1\}$ and our goal is to estimate $\Ex(Y^1 \mid V)=\gamma(V;\psi)$, where $V \subseteq L$ is a subset of the covariates. The full data orthogonal complement of the nuisance tangent space includes functions of the form
\begin{equation}
\tilde\varphi_g(Z^*;\psi) = g(V) \{ Y^1 - \gamma(V;\psi)\}
\end{equation}
for arbitrary functions $g$. From Theorem 7.2 in \cite{tsiatis2006semiparametric}, if $\pi(l)=P(A=1 \mid L=l)$ is bounded away from zero, then the observed data space $\mathcal{T}_\eta^\perp$ comprises  functions of the form
\begin{equation}
\frac{A}{\pi(L)} \Big[ \tilde\varphi_g(Z^*;\psi) +  \{1-\pi(L)\} h(Z) \Big] - (1-A) h(Z)
\end{equation}
for arbitrary functions $h$ (the simplest estimator would use the above as an estimating function with $h=0$). Note that functions of the above form only depend on observed data since $Y^1=Y$ when $A=1$.  This represents an inverse-probability-weighting approach for mapping full data spaces to observed data spaces.

\section{Empirical Processes}
\label{sec:empirical}

In the previous section we discussed how to construct influence functions $\varphi(Z;\psi,\eta)$ %(depending on observed data $Z$, parameter of  interest $\psi$, and nuisance function $\eta$) 
in semiparametric models. We also discussed how one can use these influence functions to construct estimators $\hat\psi$ for $\psi$, by solving (up to order $o_p(1/\sqrt{n})$) the estimating equation
\begin{equation}
\Pn\{\varphi(Z;\psi,\hat\eta)\}=0
\end{equation}
in $\psi$, where $\hat\eta$ is an estimator of the nuisance function. As in the previous section we let $\Pn=n^{-1} \sum_i \delta_{Z_i}$ denote the empirical measure so that sample averages can be written as $n^{-1} \sum_i f(Z_i) = \int f(z) \ \!d\Pn =\Pn\{ f(Z)\}$. We briefly discussed the asymptotics of the estimators $\hat\psi$ given above for the case where $\hat\eta \in \R^q$ is a finite-dimensional real-valued parameter, itself estimated from some estimating equation; a standard estimating equation analysis can then be used by simply stacking estimating equations for $\psi$ and $\eta$ together. 

In contrast, in this section we consider how to analyze the asymptotic behavior of $\hat\psi$ when the nuisance function $\eta$ is estimated nonparametrically, in the sense that $\hat\eta$ cannot be characterized by a finite-dimensional real-valued parameter. This can be accomplished with tools from empirical process theory. Our discussion in this section comes from work  by Andrews \cite{andrews1994empirical, andrews1994asymptotics}, Pollard \cite{pollard1984conv, pollard1990empirical}, van der Vaart \cite{vandervaart1996weak, vandervaart2000asymptotic, vandervaart2002semiparametric}, and Wellner \cite{shorack1986empirical, vandervaart1996weak}, among many others \cite{kosorok2007intro, vanderlaan2011targeted}. The field of empirical process theory is vast; we limit our discussion to tools for handling nuisance estimation.

\subsection{Motivation and Setup}
\label{subsec:motiv}

To motivate our study of empirical processes, consider our running example where the goal is to estimate the average treatment effect $\psi=\Ex(Y^1-Y^0)$. Specifically consider the doubly robust estimator for $\psi$ that solves an estimated version of the efficient influence function presented in Section \ref{subsec:eif}, i.e., the estimator given by $\hat\psi=\Pn\{m_1(Z;\hat\eta)-m_0(Z;\hat\eta)\}$ where
\begin{equation}
\label{eq:dr_est}
%\hat\psi = \Pn\left[ \frac{A \{Y - \hat\mu(L,1)\}}{\hat\pi(L)} - \frac{(1-A)\{Y-\hat\mu(L,0)\}}{1-\hat\pi(L)} + \hat\mu(L,1)-\hat\mu(L,0) \right]
m_a(Z;\eta) = m_a(Z;\pi,\mu) = \frac{I(A=a) \{Y - \mu(L,a)\}}{a\pi(L) + (1-a)\{1-\pi(L)\}} + \mu(L,a)  .
\end{equation}
Note that in this case the nuisance function is given by $\eta=(\pi,\mu)$. In observational studies the covariates $L$ are often high-dimensional, and little might be known about the propensity score and outcome regression functions $\pi$ and $\mu$, in which case it makes sense to use flexible, nonparametric, data-adaptive methods to estimate them. Of course then the asymptotic analysis presented in Section \ref{subsec:semi_influence} does not apply, since the estimators used to construct $\hat\eta=(\hat\pi,\hat\mu)$ will not be described by a single finite-dimensional parameter. Nonetheless under some conditions we can still learn about the asymptotics of $\hat\psi$ and obtain valid confidence intervals, using tools from empirical process theory.

Before going further, we need to introduce some notation. Throughout this section we will use $\Pb\{f(Z)\}= \int f(z) \ \!d\Pb$ to denote expectations of $f(Z)$ for a new observation $Z$ (treating the function $f$ as fixed); thus $\Pb\{\hat{f}(Z)\}$ is random when $\hat{f}$ is random (e.g., estimated from the sample). Contrast this with the fixed non-random quantity $\Ex\{\hat{f}(Z)\}$, which averages over randomness in both $Z$ and $\hat{f}$ and thus will not equal $\Pb\{\hat{f}(Z)\}$ except when $\hat{f}=f$ is fixed and non-random. 

Suppose for simplicity that $\hat\psi = \Pn\{m(Z;\hat\eta)\}$ for some $m$, as in the above example. If we only have $\Pn\{\varphi(Z;\hat\psi,\hat\eta)\}=0$ then we can proceed similarly, with an extra step requiring differentiability of $\Pb\{\varphi(Z;\psi,\eta)\}$ in $\psi$, at $\psi_0$ in a neighborhood of $\eta_0$ \cite{vandervaart2000asymptotic}. Also suppose that $\Pb\{m(Z;\eta_0)\}=\psi_0$ (alternatively we can define $\psi_0$ so that this holds by definition). For instance, it is straightforward to check for the doubly robust estimator described above that $\Pb\{m(Z;\pi_0,\mu)\}=\Pb\{m(Z;\pi,\mu_0)\}=\psi_0$ where $m=m_1-m_0$. Then consider the decomposition
\begin{align}
\label{eq:decomp}
\hat\psi - \psi_0 &= \Pn\{m(Z;\hat\eta)\} - \Pb\{m(Z;\eta_0)\} \\
&= (\Pn-\Pb)m(Z;\hat\eta) + \Pb\{m(Z;\hat\eta) - m(Z;\eta_0)\}, \nonumber
\end{align}
where the first line is true by definition, and the second follows by simply adding and subtracting $\Pb\{m(Z;\hat\eta)\}$. 

We will show that the first term $(\Pn-\Pb)m(Z;\hat\eta)$ above can be handled under general conditions with empirical process theory. Specifically, we will discuss conditions under which 
\begin{equation}
\label{eq:term1}
(\Pn-\Pb)m(Z;\hat\eta) = (\Pn-\Pb)m(Z;\eta_0) + o_p(1/\sqrt{n}),
\end{equation}
where $\hat\eta$ converges to $\eta_0$, so that $(\Pn-\Pb)m(Z;\hat\eta)$ is asymptotically equivalent to its limiting version $(\Pn-\Pb)m(Z;\eta_0)$ (up to order $o_p(1/\sqrt{n})$) and can be analyzed with a standard central limit theorem. The second term in the decomposition in \eqref{eq:decomp} typically requires a case-by-case analysis, but we will give examples shortly. Note that if we have $\Pb\{m(Z;\hat\eta) - m(Z;\eta_0)\}= (\Pn-\Pb) \phi(Z;\eta_0) + o_p(1/\sqrt{n})$ for some finite-variance function $\phi$, then
\begin{equation}
\hat\psi - \psi_0 = (\Pn-\Pb)\{m(Z;\eta_0) + \phi(Z;\eta_0)\} + o_p(1/\sqrt{n})
\end{equation}
and thus $\hat\psi$ is regular and asymptotically linear with influence function $(m+\phi)$.

\subsection{Donsker Classes}
\label{subsec:donsker}

From an empirical process perspective, a primary way to control how close the term $(\Pn-\Pb)m(Z;\hat\eta)$ is to its limiting version $(\Pn-\Pb)m(Z;\eta_0)$ (in large samples) is to restrict the complexity of the nuisance function $\eta_0$ and its estimator $\hat\eta$. If these functions are not too complex, then the terms will not differ by more than $o_p(1/\sqrt{n})$. In this subsection we will discuss characterizing complexity with Donsker classes.

We will start by giving the main result in the context of our example, and will then describe the conditions in detail. Suppose our nuisance estimator $\hat\eta$ converges to some limit $\eta_0$ in the sense that 
\begin{equation}
|| m(;\hat\eta)-m(;\eta_0) ||^2 = \int \{ m(z;\hat\eta)-m(z;\eta_0)\}^2 \ \! dP(z) = o_p(1) ,
\end{equation}
and suppose the function class  $\mathcal{M}=\{ m(; \eta) :  \eta \in H \}$ is a {Donsker class} (to be defined shortly), where $H$ is a function class containing the nuisance estimator $\hat\eta$. Then the result in \eqref{eq:term1} holds, i.e.,
\begin{equation}
(\Pn-\Pb)m(Z;\hat\eta) = (\Pn-\Pb)m(Z;\eta_0) + o_p(1/\sqrt{n}).
\end{equation}
Thus, asymptotically, nuisance estimation only affects the second term in \eqref{eq:decomp}.

In order to define a Donsker class, we need to introduce a few concepts first. Throughout this section we use $\Gn=\sqrt{n}(\Pn-\Pb)$ for ease of notation. Let $\mathcal{F}$ denote a class of functions $f: \mathcal{Z} \rightarrow \R$, and consider the \textit{empirical process}
\begin{equation}
\{ \Gn f : f \in \mathcal{F} \} .
\end{equation}
This is a type of \textit{stochastic process} since it is a collection of random variables indexed by a set (the function class $\mathcal{F}$). From one standpoint, given a function $f$, we can view $\Gn f=\sqrt{n}(\Pn-\Pb)f(Z)$ as a random variable mapping the sample (product) space $\mathcal{Z}^n$ to $\R$. Alternatively, given a sample $(Z_1,...,Z_n)$, we can also view $\Gn f$ as a map from the function class $\mathcal{F}$ to $\R$. Therefore (if these latter maps are bounded) we can view the empirical process as a \textit{random function}, mapping the sample space $\mathcal{Z}^n$ to the space $\ell^\infty(\mathcal{F})$ of bounded functions $h: \mathcal{F} \rightarrow \R$ with $\sup_{f \in \mathcal{F}} |h(f)| = ||h||_{\mathcal{F}} < \infty$. 

The above discussion of the empirical process $\{\Gn f : f \in \mathcal{F}\}$ was all for a fixed sample size $n$. Now consider a sequence of empirical processes $\{\Gn f : f \in \mathcal{F}\}_{n \geq 1}$. We say this sequence \textit{converges in distribution} to element $\Gb$ (equivalently, converges weakly to $\Gb$) in the space $\ell^\infty(\mathcal{F})$, denoted $\Gn \indist \Gb$, if 
\begin{equation}
\Ex^* h(\Gn) \rightarrow \Ex h(\Gb) 
\end{equation}
for all continuous bounded functions $h: \ell^\infty(\mathcal{F}) \rightarrow \R$, where $\Ex^*$ denotes outer expectation. (Outer expectation is a measure-theoretic subtlety that we will largely sidestep here; roughly, $\Ex^*$ can be viewed as a generalization of expectation that accounts for the fact that $h(\Gn)$ may not be measurable). Thus we have a notion of convergence for empirical processes viewed as random functions. Finally, we say a generic measurable random element $\Gb$ is \textit{tight} if for all $\epsilon>0$ there is a compact set $S$ for which $P(\Gb \in S)>1-\epsilon$, i.e., if the element $\Gb$ stays in a compact set with high probability.

We are now ready to define a {Donsker class}. A function class $\mathcal{F}$ is called a \textit{Donsker class} if the sequence of empirical processes $\{\Gn f : f \in \mathcal{F}\}_{n \geq 1}$ converges in distribution to some tight limit $\Gb$ (in fact this limit must be a zero-mean Gaussian process $\Gb_P$, known as a $P$-Brownian bridge).

The Donsker property, along with the continuous mapping theorem, allow us to obtain results like that given in \eqref{eq:term1}. Specifically, suppose $\hat{f} \in \mathcal{F}$ for a Donsker class $\mathcal{F}$, and suppose $\hat{f}$ converges to $f_0$ in the sese that $||\hat{f}-f_0||=o_p(1)$, where $||f||^2=Pf^2$ denotes the $L_2(P)$ norm as before. Then (as in Lemma 19.24 of \cite{vandervaart2000asymptotic}) we can apply the continuous mapping theorem to $(\Gn, \hat{f}) \indist (\Gb_P,f_0)$ with function $h(z,f)=z(f)-z(f_0)$ to obtain that
\begin{equation}
\label{eq:lemma}
\Gn \hat{f} = \Gn f_0 + o_p(1) .
\end{equation}
Thus $(\Pn-\Pb) \hat{f} = n^{-1/2} \Gn \hat{f}$ is asymptotically equivalent to $ (\Pn-\Pb) f_0$, up to $o_p(1/\sqrt{n})$ error.

In our setting, where $\hat\psi=\Pn\{m(Z;\hat\eta)\}$, it is often more natural to put Donsker conditions on the estimated nuisance functions themselves, i.e., to assume that $\hat\eta \in H$ for a Donsker class $H$, rather than to put conditions on the transformed function class $\mathcal{M} =\{m(;\eta) : \eta \in H\}$. Fortunately, `nice enough' transformations of Donsker function classes will also be Donsker. Specifically, suppose the function classes $\mathcal{F}$ and $\mathcal{F}_j$ are Donsker; then, as discussed in Section 2.10 of \cite{vandervaart1996weak}, as in \cite{andrews1994empirical, vandervaart2000asymptotic}, the following transformations of $\mathcal{F}$ and $\mathcal{F}_j$ are also Donsker:
\begin{enumerate}
\item \textit{Subsets:} $\mathcal{G} \subset \mathcal{F}$
\item \textit{Unions:} $\mathcal{G} = \mathcal{F}_1 \cup \mathcal{F}_2$
\item \textit{Closures:} $\mathcal{G} =\{ g : f_m \rightarrow g \ \text{ pointwise and in $L_2$, for } \ f_m \in \mathcal{F} \}$
\item \textit{Convex combinations:} $\mathcal{G} = \{ g : g = \sum_i w_i  f_i \ \text{ for } \ f_i \in \mathcal{F}, \sum_i |w_i| \leq 1 \}$
\item \textit{Lipschitz transformations:} $\mathcal{G}=\{ g : g=\phi(f_1,...,f_k) \ \text{ for } \ f_j \in \mathcal{F}_j \}$ if $\phi$ satisfies $|\phi(f_1,...,f_k)(x)-\phi(f'_1,...,f'_k)(x)|^2 \leq \sum_j (f_j-f'_j)(x)^2$ for all $f_j$, $f'_j$, and $x$, and if $\sup_{f \in \mathcal{F}_j} |P f| < \infty$ and $\int \phi(f_1,...,f_k)(x)^2 dx<\infty$.
\end{enumerate}

The convex combination result suggests using ensemble methods that use weighted combinations of estimators, e.g., Super Learner \cite{vanderlaan2003unified2, vanderlaan2007super, vanderlaan2011targeted}. The Lipschitz transformation result given above is particularly useful. It means, for example, that the following function classes are Donsker \cite{andrews1994empirical, vandervaart1996weak, vandervaart2000asymptotic}:
\begin{enumerate}
\item \textit{Minimums:} $\mathcal{G} = \{ g : g=\min(f_1,f_2) \ \text{ for } \ f_j \in \mathcal{F}_j \}$
\item \textit{Maximums:} $\mathcal{G} = \{ g : g=\max(f_1,f_2) \ \text{ for } \ f_j \in \mathcal{F}_j \}$
\item \textit{Sums:} $\mathcal{G} = \{ g : g= f_1+f_2 \ \text{ for } \ f_j \in \mathcal{F}_j \}$
\item \textit{Products:} $\mathcal{G} = \{ g : g = f_1 f_2 \ \text{ for } \ f_j \in \mathcal{F}_j\}$ if $\mathcal{F}_j$ are uniformly bounded
\item \textit{Ratios:} $\mathcal{G} = \{ g: g= 1/f \ \text{ for } \ f \in \mathcal{F}\}$ if $f \geq \delta>0$ for all $f \in \mathcal{F}$
\end{enumerate}
Repeated use of stability results like those above often allows one to conclude Donsker properties for the class $\mathcal{M} =\{m(;\eta) : \eta \in H\}$ based on Donsker assumptions about the class $H$. 

For example, consider the doubly robust estimator $\hat\psi = \Pn\{m_1(Z;\hat\eta)-m_0(Z;\hat\eta)\}$ given in \eqref{eq:dr_est}. If $\hat\pi$ and $\hat\mu$ take values in Donsker classes $\mathcal{F}_\pi$ and $\mathcal{F}_\mu$, respectively, then $m_a(Z;\hat\eta)$ does as well (provided that $\pi$ is bounded away from zero and one for all $\pi \in \mathcal{F}_\pi$). This follows from Lipschitz results 3 and 5 for sums and ratios above.

\subsection{Examples of Donsker Classes}
\label{subsec:donsker_ex}

To this point we have seen that, if we assume the estimated nuisance functions $\hat\eta$ are contained in Donsker function classes, we can use a standard central limit theorem to analyze $(\Pn-\Pb) m(Z;\hat\eta)$ since it is asymptotically equivalent to $(\Pn-\Pb) m(Z;\eta_0)$ up to order $o_p(1/\sqrt{n})$. We have defined Donsker classes and shown how they can be combined and modified to produce new Donsker classes, but we have yet to give any specific examples of such classes. For the prior results to be useful over and above more standard parametric techniques, we need Donsker classes to be able to capture sufficiently flexible functions. Luckily, this is in fact the case, as we will discuss in this subsection using specific examples.

First we will simply provide a short list of function classes that are Donsker, and then we will briefly discuss how one typically shows that a particular class is Donsker (using bracketing and covering numbers). Results showing that certain classes are Donsker are somewhat scattered across the literature, but examples and nice overviews are given by \cite{vandervaart1996weak,vandervaart2000asymptotic}, for example. Among many other kinds of classes, the following simple classes of functions are Donsker classes \cite{gill1995inefficient, vandervaart1996weak,vandervaart2000asymptotic}:
\begin{enumerate}
\item \textit{Indicator functions:} $\mathcal{F}=\{ f: f(x) = I(x<t) , t \in \R \}$
\item \textit{Vapnik-Cervonenkis (VC) classes}
\item \textit{Bounded monotone functions}
\item \textit{Lipschitz parametric functions:} $\mathcal{F}=\{ f: f(x) = f(x;\theta) , \theta \in \Theta \subset \R^q \}$ with $|f(x;\theta_1)-f(x;\theta_2)| \leq b(x) ||\theta_1-\theta_2||$ for some $b$ with $\int |b(x)|^r \ \! dP(x) < \infty$
\item \textit{Smooth functions:} $\mathcal{F}=\{ f: \sup_x |\frac{\partial^\alpha f(x_1,...,x_q)}{\partial^{\alpha_1} x_1 ... \partial^{\alpha_q} x_q}| < B < \infty, \ \text{ with } \ \alpha>q/2 \}$
\item \textit{Sobolev classes:} $\{ f : \sup_x |f(x)| \leq 1 , f^{(k-1)} \text{ absolutely cts.}, \int |f^{(k)}(x)|^2 \ dx \leq 1\}$
\item \textit{Uniform sectional variation:} $\{ f : \sup_{x_1} ||f(x_1,\cdot)||_{tv} \leq B_1, \sup_{x_2} ||f(\cdot,x_2)||_{tv} \leq B_2\}$ where $B_1, B_2 < \infty$ and $||\cdot||_{tv}$ denotes the total variation norm.
\end{enumerate}

Thus we see that Donsker classes include usual parametric classes, but many other classes as well, including infinite-dimensional classes that only require certain smoothness or boundedness. Many other function classes can also be shown to be Donsker. For example, any appropriate combination or transformation of the above classes as discussed in the previous subsection will also be Donsker. 

Showing that a function class is Donsker is often accomplished using bracketing or covering numbers \cite{vandervaart1996weak, vandervaart2000asymptotic}, which are measures of the size of a class $\mathcal{F}$. These measures also provide simple sufficient conditions for a function class being Donsker. An $\epsilon$-bracket (in $L_2(P)$) is defined as all functions $f$  bracketed by functions $[l,u]$ (i.e., $l \leq f \leq u$) satisfying $\int \{u(z)-l(z)\}^2 \ \!dP(z) < \epsilon^2$.  The \textit{bracketing number} of a class $\mathcal{F}$ is the smallest number of $\epsilon$-brackets needed to cover $\mathcal{F}$, and is denoted by $N_B(\epsilon,\mathcal{F})$. Similarly, the \textit{covering number} of a class $\mathcal{F}$ (with envelope $F$, i.e., $\sup_{\mathcal{F}} |f| \leq F$) is the smallest number of $L_2(Q)$ balls of radius $\epsilon$ needed to cover $\mathcal{F}$, and is denoted by $N_C(\epsilon,\mathcal{F})$. Then the class $\mathcal{F}$ is Donsker if either 
\begin{equation}
\int_0^1 \sqrt{\log N_B(\epsilon,\mathcal{F}) } \ d\epsilon < \infty, \ \ \text{  or } \ \int_0^1 \sqrt{ \log \sup_Q N_C(\epsilon \sqrt{QF^2},\mathcal{F}) } \ d\epsilon < \infty .
\end{equation}

\subsection{Average Treatment Effect Example}
\label{subsec:ate_ex}

Now we return to analyze the asymptotic behavior of the doubly robust estimator of the average treatment effect $\psi=\E(Y^1-Y^0)$ from Section \ref{subsec:eif}, which is given by $\hat\psi=\Pn\{m(Z;\hat\eta)\}= \Pn\{m_1(Z;\hat\eta)-m_0(Z;\hat\eta)\}$ with
\begin{equation}
\label{eq:dr_est}
%\hat\psi = \Pn\left[ \frac{A \{Y - \hat\mu(L,1)\}}{\hat\pi(L)} - \frac{(1-A)\{Y-\hat\mu(L,0)\}}{1-\hat\pi(L)} + \hat\mu(L,1)-\hat\mu(L,0) \right]
m_a(Z;\eta) = m_a(Z;\pi,\mu) = \frac{I(A=a) \{Y - \mu(L,a)\}}{a\pi(L) + (1-a)\{1-\pi(L)\}} + \mu(L,a)  .
\end{equation}

Throughout we assume the identification assumptions from Section \ref{subsec:semi_id}, or else suppose we are estimating the observed data quantity $\Ex\{\mu(L,1)-\mu(L,0)\}$ under the positivity assumption. Suppose the estimator $\hat\eta=(\hat\pi,\hat\mu)$ converges to some $\overline\eta=(\overline\pi,\overline\mu)$ in the sense that $||\hat\eta-\overline\eta||=o_p(1)$, where either $\overline\pi=\pi_0$ or $\overline\mu=\mu_0$ (or both) correspond to the true nuisance function. Thus at least one nuisance estimator needs to converge to the correct function, but one can be misspecified. Then $\Pb\{m(Z;\overline\eta)\}=\Pb\{m(Z;\eta_0)\}=\psi_0$, from the easy-to-check fact that $\Pb\{m(Z;\pi_0,\mu)\}=\Pb\{m(Z;\pi,\mu_0)\}$ for any $\overline\pi$ and $\overline\mu$. Thus as in Section \ref{subsec:motiv} we can write
\begin{equation}
\hat\psi - \psi_0  = (\Pn-\Pb)m(Z;\hat\eta) + \Pb\{m(Z;\hat\eta) - m(Z;\overline\eta)\} .
\end{equation}
As discussed in Section \ref{subsec:donsker}, if the estimators $\hat\pi$ and $\hat\mu$ take values in Donsker classes, then $m_a(Z;\hat\eta)$ does as well (as long as functions in the class containing $\hat\pi$ are uniformly bounded away from zero and one). Therefore the result in \eqref{eq:term1} applies, and we have
\begin{equation}
\hat\psi - \psi_0  = (\Pn-\Pb)m(Z;\overline\eta) + \Pb\{m(Z;\hat\eta) - m(Z;\overline\eta)\} + o_p(1/\sqrt{n}) .
\end{equation}

Now it remains to analyze $\Pb\{m(Z;\hat\eta) - m(Z;\overline\eta)\}$. By iterated expectation this term equals
\begin{equation}
\sum_{a \in \{0,1\}} \Pb \left[ \frac{\pi_0(L)-\hat\pi(L)}{a \hat\pi(L) + (1-a) \{1-\hat\pi(L)\}} \{\mu_0(L,a)-\hat\mu(L,a)\} \right] .
\end{equation}
Therefore, by the fact that $\hat\pi$ is bounded away from zero and one, along with the Cauchy-Schwarz inequality ($P(fg) \leq ||f|| \ ||g||$), we have that (up to a multiplicative constant) $|\Pb\{m(Z;\hat\eta) - m(Z;\overline\eta)\}|$ is bounded above by
\begin{equation}
\sum_{a \in \{0,1\}} || \pi_0(L)-\hat\pi(L)|| \ || \mu_0(L,a)-\hat\mu(L,a) || .
\end{equation}
Thus for example if $\hat\pi$ is based on a correctly specified parametric model, so that $||\hat\pi-\pi_0||=O_p(n^{-1/2})$, then we only need $\hat\mu$ to be consistent, $||\hat\mu-\mu_0||=o_p(1)$, to make the product term $\Pb\{m(Z;\hat\eta) - m(Z;\overline\eta)\}=o_p(1/\sqrt{n})$ asymptotically negligible. Then the doubly robust estimator satisfies $\hat\psi-\psi_0=(\Pn-\Pb)m(Z;\eta_0) + o_p(1/\sqrt{n})$ and it is efficient with influence function $\varphi(Z;\psi,\eta)=m(Z;\eta)-\psi$. Thus if we know the treatment mechanism, the outcome model can be very flexible.

Another way to achieve efficiency is if we have both $||\hat\pi-\pi_0||=o_p(n^{-1/4})$ and $||\hat\mu-\mu_0||=o_p(n^{-1/4})$, so that the product term is $o_p(1/\sqrt{n})$ and asymptotically negligible. This of course occurs if both $\hat\pi$ and $\hat\mu$ are based on correctly specified models, but it can also hold even for estimators that are very flexible and not based on parametric models. However, completely nonparametric (e.g., kernel or nearest-neighbor) estimators are typically not an option in this setting since they will generally converge at rates slower than $n^{-1/4}$; exceptions include cases where there are very few covariates or very strong smoothness assumptions. Explicit conditions ensuring given convergence rates for kernel estimators are described for example in \cite{newey1994large}. Thus some modeling is in general required to attain $n^{-1/4}$ rates, but luckily numerous semiparametric models yield estimators that can satisfy this condition. In particular, faster than $n^{-1/4}$ rates are possible with single index models, generalized additive models, and partially linear models (see for example \cite{horowitz2009semiparametric} for a review of such models, which typically yield estimators with $n^{-2/5}$ rates), as well as regularized estimators such as the Lasso  \cite{belloni2014inference, belloni2015uniform}. Cross-validation-based weighted combinations of such estimators (e.g., Super Learner) can also satisfy this rate condition if one of the candidate estimators does \cite{vanderlaan2003unified2}. 

Inference after nonparametric estimation of $\eta$ in truly doubly robust settings where one arbitrary nuisance estimator can be misspecified is more complicated. If one of the estimators $\hat\pi$ or $\hat\mu$ is misspecified so that either $||\hat\pi-\pi_0||=O_p(1)$ or $||\hat\mu-\mu_0||=O_p(1)$, then obtaining root-n rate inference for standard estimators will typically require knowledge of which estimator is correctly specified, as well as that the correctly specified estimator is based on a parametric model. More sophisticated estimators that weaken this requirement are discussed in the next section (e.g., \cite{vanderlaan2014targeted2}).

\section{Extensions \& Future Directions}
\label{sec:future}

In this section we briefly describe some future directions and extensions to semiparametric causal inference beyond the theory we have presented in this review. A number of authors have worked to extend semiparametric causal inference to, for example, settings involving non-standard sampling, estimation and inference under yet weaker conditions on the nuisance estimators, and complex non-regular or non-smooth parameters. 

Throughout this review we presumed access to an independent and identically distributed sample from the distribution $P$ of interest; however, many studies use alternative sampling schemes. For example, authors have developed results for semiparametric causal inference in case control studies \cite{vanderlaan2008estimation, vanderweele2011weighting, tchetgen2012double,  rose2014double, vanderweele2014invited} and matched cohort studies \cite{vanderlaan2013estimating, kennedy2015semiparametric}. There has also been progress made for causal inference in studies using network data with possible interference \cite{hudgens2012towards, tchetgen2012causal, ogburn2014causal, vanderlaan2014causal}. Much more work is needed in settings related to both study designs with non-standard sampling and network data with interference. The latter should be a growing concern as data from, e.g., social networks becomes more commonplace. 

In Section \ref{sec:empirical} we showed that semiparametric estimators can have appealing asymptotic behavior, including standard root-n rates of convergence and straightforward confidence intervals, even when using flexible nonparametric estimates of nuisance functions. However, as noted in Section \ref{subsec:ate_ex}, this can require a delicate balance in settings where one does not want to rely on parametric models, and also wants to be agnostic about whether the treatment or outcome process is correctly estimated. Efforts to weaken the conditions needed on the nuisance estimation have been made using approaches based on higher-order estimation \cite{vanderlaan2014targeted2, carone2014higher, diaz2015second}, which were inspired by work by Robins et al. \cite{robins2008higher, robins2009quadratic, vandervaart2014higher} that focused on minimax estimation in settings where root-n rates of convergence are not possible. Further, Donsker-type regularity conditions (though not rate conditions) can be weakened via cross-validation approaches, proposed for example by \cite{zheng2010asymptotic}.

We also supposed in this review that our target parameter was a low-dimensional Euclidean parameter $\psi \in \R^p$ that admitted regular asymptotically linear estimators. However, in some settings these conditions fail to hold. As mentioned above, Robins et al. \cite{robins2008higher, robins2009quadratic, vandervaart2014higher} considered semiparametric minimax estimation in settings where the parameter of interest is Euclidean, but root-n rates of convergence cannot be attained due to high-dimensional covariates. Estimation of functional effect parameters was considered by \cite{diaz2013targeted, kennedy2015nonparametric} in the context of continuous treatment effects; in such settings the target parameter is a non-pathwise differentiable curve, and root-n rates of convergence are again not possible. Inference for a non-regular parameter in an optimal treatment regime setting was considered by \cite{luedtke2014statistical}; in this case non-regularity does not preclude the existence of root-n rate inference.

Numerous other authors have also made important contributions extending semiparametric causal inference to novel settings; unfortunately we cannot list all of them here. In addition, much important work is left to be done, both in the areas mentioned above as well as in many other interesting settings.

\begin{acknowledgement}
Edward Kennedy acknowledges support from NIH grant R01-DK090385, and thanks Jason Roy and Bret Zeldow for very helpful comments and discussion.
\end{acknowledgement}


\begin{thebibliography}{99.}%
% Use the following syntax and markup for your references if the subject of your book is from the field "Mathematics, Physics, Statistics, Computer Science"
% Contribution 
%\bibitem{science-contrib} Broy, M.: Software engineering --- from auxiliary to key technologies. In: Broy, M., Dener, E. (eds.) Software Pioneers, pp. 10-13. Springer, Heidelberg (2002)
% Online Document
%\bibitem{science-online} Dod, J.: Effective substances. In: The Dictionary of Substances and Their Effects. Royal Society of Chemistry (1999) Available via DIALOG. \\ \url{http://www.rsc.org/dose/title of subordinate document. Cited 15 Jan 1999}
% Monograph
%\bibitem{science-mono} Geddes, K.O., Czapor, S.R., Labahn, G.: Algorithms for Computer Algebra. Kluwer, Boston (1992) 
% Journal article
%\bibitem{science-journal} Hamburger, C.: Quasimonotonicity, regularity and duality for nonlinear systems of partial differential equations. Ann. Mat. Pura. Appl. \textbf{169}, 321--354 (1995)
% Journal article by DOI
%\bibitem{science-DOI} Slifka, M.K., Whitton, J.L.: Clinical implications of dysregulated cytokine production. J. Mol. Med. (2000) doi: 10.1007/s001090000086 

\bibitem{andrews1994empirical} Andrews, D.W.K.: Empirical process methods in econometrics. Handbook of Econometrics \textbf{4}, 2247--2294 (1994)

\bibitem{andrews1994asymptotics}Andrews, D.W.K.: Asymptotics for semiparametric econometric models via stochastic equicontinuity. Econometrica, 43--72 (1994)

\bibitem{angrist1996ident} Angrist, J.D., Imbens, G.W., Rubin, D.B.: Identification of causal effects using instrumental variables. J. Am. Stat. Assoc. \textbf{91}, 444--455 (1996)

\bibitem{begun1983information} Begun, J.M., Hall, W.J., Huang, W.M., Wellner, J.A.: Information and asymptotic efficiency in parametric-nonparametric models. Ann. Stat. \textbf{11}, 432--452 (1983)

\bibitem{belloni2014inference} Belloni, A., Chernozhukov, V., Hansen, C.: Inference on treatment effects after selection among high-dimensional controls. Rev. Econ. Stud. \textbf{81}, 608--650 (2014)

\bibitem{belloni2015uniform} Belloni, A., Chernozhukov, V., Kato, K.: Uniform post-selection inference for least absolute deviation regression and other Z-estimation problems. Biometrika. \textbf{102}, 77-94 (2015)

\bibitem{bickel1993efficient} Bickel, P.J., Klaassen, C.A.J., Ritov, Y., Wellner, J.A.: Efficient and Adaptive Estimation for Semiparametric Models. Springer, New York (1993) 

\bibitem{boos2002calculus} Stefanski, L.A., Boos, D.D.: The calculus of M-estimation. Am. Stat. \textbf{56}, 29--38 (2002)

\bibitem{chakraborty2013statistical} Chakraborty, B., Moodie, E.E.M.: Statistical Methods for Dynamic Treatment Regimes. Springer, New York. (2013)

\bibitem{carone2014higher} Carone, M., Diaz, I., van der Laan, M.J.: Higher-order targeted minimum loss-based estimation. U.C. Berkeley Division of Biostatistics Working Paper Series. \textbf{331}, 1--39 (2015)

\bibitem{dawid2000causal} Dawid, P.A.: Causal inference without counterfactuals. J. Am. Stat. Assoc. \textbf{95}, 407--424 (2000)

\bibitem{diaz2013targeted} Diaz, I., van der Laan, M.J.: Targeted data adaptive estimation of the causal dose–response curve. J. Causal Inf. \textbf{1}, 171--192 (2013)

\bibitem{diaz2015second} Diaz, I., Carone, M., van der Laan, M.J.: Second order inference for the mean of a variable missing at random. U.C. Berkeley Division of Biostatistics Working Paper Series. \textbf{337}, 1--22 (2015)

\bibitem{gill1995inefficient} Gill, R.D., van der Laan, M.J., Wellner, J.A. Inefficient estimators of the bivariate survival function for three models. Ann. Inst. Henri Poincare. \textbf{31},  545--597. (1995)

\bibitem{gill1997car} Gill, R.D., van der Laan, M.J., Robins, J.M. Coarsening at random: Characterizations, conjectures, counter-examples. In: Proceedings of the First Seattle Symposium in Biostatistics, pp. 255--294. Springer, New York (1997)

\bibitem{hahn1998role} Hahn, J.: On the role of the propensity score in efficient semiparametric estimation of average treatment effects. Econometrica. \textbf{66}, 315--333. (1998)

\bibitem{hernan2006instruments} Hernan, M.A., Robins, J.M.: Instruments for causal inference: an epidemiologist's dream?. Epidemiology. \textbf{17}, 360--372 (2006)

\bibitem{horowitz2009semiparametric} Horowitz, J.L.: Semiparametric and Nonparametric Methods in Econometrics. Springer, New York (2009)

\bibitem{hudgens2012towards} Hudgens, M.G., Halloran, M.E.: Toward causal inference with interference. J. Am. Stat. Assoc. \textbf{103}, 832--842  (2012)

\bibitem{kennedy2015semiparametric} Kennedy, E.H., Sjolander, A., Small, D.S.: Semiparametric causal inference in matched cohort studies. Biometrika. \textbf{102}, 739-746 (2015)

\bibitem{kennedy2015nonparametric} Kennedy, E.H., Ma, Z., McHugh, M.D., Small, D.S.: Nonparametric methods for doubly robust estimation of continuous treatment effects. arXiv preprint, arXiv:1507.00747 (2015)

\bibitem{kosorok2007intro} Kosorok, M.R.: Introduction to Empirical Processes and Semiparametric Inference. Springer, New York (2007) 

\bibitem{luedtke2014statistical} Luedtke, A.R., van der Laan, M.J.: Statistical inference for the mean outcome under a possibly non-unique optimal treatment strategy. U.C. Berkeley Division of Biostatistics Working Paper Series. \textbf{332}, 1--37 (2014)

\bibitem{manski2003partial} Manski, C.F.: Partial Identification of Probability Distributions. Springer, New York (2003)

\bibitem{murphy2003optimal} Murphy, S.A.: Optimal dynamic treatment regimes. J. Roy. Stat. Soc. B. \textbf{65}, 331--355 (2003)

\bibitem{neugebauer2007nonparametric} Neugebauer, R., van der Laan, M.J.: Nonparametric causal effects based on marginal structural models. J. Stat. Plan. Infer. \textbf{137}, 419--434 (2007)

\bibitem{newey1994large} Newey, W.K., McFadden, D.: (1994). Large sample estimation and hypothesis testing. Handbook of Econometrics \textbf{4}, 2111--2245 (1994)

\bibitem{newey1994asymptotic} Newey, W.K.: The asymptotic variance of semiparametric estimators. Econometrica. \textbf{62}, 1349--1382. (1994)

\bibitem{neyman1923sur} Neyman, J.: On the application of probability theory to agricultural experiments: Essay on principles. Excerpts reprinted (1990) in English (D. Dabrowska and T. Speed, trans.). Statist. Sci. \textbf{5}, 463--472 (1923)

\bibitem{ogburn2014causal} Ogburn, E.L., VanderWeele, T.J.: Causal diagrams for interference. Stat. Sci. \textbf{29}, 559--578 (2014)

\bibitem{pearl1995causal} Pearl, J.: Causal diagrams for empirical research. Biometrika \textbf{82}, 669--688 (1995)

\bibitem{pearl2009causality} Pearl, J.: Causality. Cambridge University Press (2009)

\bibitem{petersen2010diagnosing} Petersen, M.L., Porter, K.E., Gruber, S., Wang, Y., van der Laan, M.J.: Diagnosing and responding to violations in the positivity assumption. Stat. Methods. Med. Res. \textbf{21}, 31--54 (2010)

\bibitem{pfanzagl1982contributions} Pfanzagl, J.: Contributions to a General Asymptotic Statistical Theory. Springer, New York (1982) 

\bibitem{pfanzagl1990estimation} Pfanzagl, J.: Estimation in Semiparametric Models. Springer, New York (1990)

\bibitem{pollard1984conv} Pollard, D.: Convergence of stochastic processes. Springer, New York (1984) 

\bibitem{pollard1990empirical} Pollard, D.: Empirical processes: theory and applications. In NSF-CBMS Regional Conference Series in Probability and Statistics. Institute of Mathematical Statistics and the American Statistical Association (1990)

\bibitem{robins1986new} Robins, J.M.: A new approach to causal inference in mortality studies with a sustained exposure period -- application to control of the healthy worker survivor effect. Math. Mod. \textbf{7}, 1393--1512 (1986)

\bibitem{robins1994estimation} Robins, J.M., Rotnitzky, A., Zhao, L.P.: Estimation of regression coefficients when some regressors are not always observed. J. Am. Stat. Assoc. \textbf{89}, 846--866 (1994)

\bibitem{robins1995analysis} Robins, J.M., Rotnitzky, A., Zhao, L.P.: Analysis of semiparametric regression models for repeated outcomes in the presence of missing data. J. Am. Stat. Assoc. \textbf{90}, 106--121 (1995)

\bibitem{robins1999sensitivity} Robins, J.M., Rotnitzky, A., Scharfstein, D.O.: Sensitivity analysis for selection bias and unmeasured confounding in missing data and causal inference models. In: Statistical Models in Epidemiology, the Environment, and Clinical Trials, pp. 1--94. Springer, New York (1999)

\bibitem{robins2000marginal} Robins, J.M., Hernan, M.A., Brumback, B.: Marginal structural models and causal inference in epidemiology. Epidemiology. \textbf{11}, 550--560 (2000)

\bibitem{robins2008higher} Robins, J.M., Li, L., Tchetgen Tchetgen, E.J., van der Vaart, A.W.: Higher order influence functions and minimax estimation of nonlinear functionals. In: Probability and Statistics: Essays in Honor of David A. Freedman, pp. 335--421. Institute of Mathematical Statistics (2008)

\bibitem{robins2009estimation} Robins, J.M., Hernan, M.A.: Estimation of the causal effects of time-varying exposures. In: Fitzmaurice, G., Davidian, M., Verbeke, G., Molenberghs, G. (eds.) Longitudinal Data Analysis, pp. 553--600. Chapman \& Hall, London (2009)

\bibitem{robins2009quadratic} Robins, J.M., Li, L., Tchetgen Tchetgen, E.J., van der Vaart, A.W.: Quadratic semiparametric von mises calculus. Metrika. \textbf{69}, 227--247 (2009)

\bibitem{rose2014double} Rose, S., van der Laan, M.J.: A double robust approach to causal effects in case-control studies. Am. J. Epid. \textbf{179}, 662--669  (2014)

\bibitem{rubin1974estimating} Rubin, D.B.: Estimating causal effects of treatments in randomized and nonrandomized studies. J. Educ. Psychol. \textbf{66}, 688--701 (1974)

\bibitem{rubin1978role} Rubin, D.B.: Bayesian inference for causal effects: The role of randomization. Ann. Stat. \textbf{6}, 34--58 (1978)

\bibitem{shorack1986empirical} Shorack, G.R., Wellner, J.A.: Empirical Processes with Applications to Statistics. Wiley, New York (1986)

\bibitem{tchetgen2012double} Tchetgen Tchetgen, E.J., Rotnitzky, A.: Double-robust estimation of an exposure-outcome odds ratio adjusting for confounding in cohort and case-control studies. Stat. Med. \textbf{30}, 335--347 (2011)

\bibitem{tchetgen2012semiparametric} Tchetgen Tchetgen, E.J., Shpitser, I.: Semiparametric theory for causal mediation analysis: efficiency bounds, multiple robustness and sensitivity analysis. Ann. Stat. \textbf{40}, 1816--1845 (2012)

\bibitem{tchetgen2012causal} Tchetgen Tchetgen, E.J., VanderWeele, T. J.: On causal inference in the presence of interference. Stat. Methods Med. Res. \textbf{21}, 55--75 (2012)

\bibitem{tsiatis2006semiparametric} Tsiatis, A.A.: Semiparametric Theory and Missing Data. Springer, New York (2006) 

\bibitem{vanderlaan2003unified} van der Laan, M.J., Robins, J.M.: Unified Methods for Censored Longitudinal Data and Causality. Springer, New York (2003) 

\bibitem{vanderlaan2003unified2} van der Laan, M.J., Dudoit, S.: Unified cross-validation methodology for selection among estimators and a general cross-validated adaptive epsilon-net estimator: Finite sample oracle inequalities and examples. U.C. Berkeley Division of Biostatistics Working Paper Series. \textbf{130}, 1--103 (2003)

\bibitem{vanderlaan2006targeted} van der Laan, M.J., Rubin, D.: Targeted maximum likelihood learning. Int. J. Biostat. \textbf{2}, 1--38 (2006)

\bibitem{vanderlaan2007super} van der Laan, M.J., Polley, E.C., Hubbard, A.E.: Super learner. Stat. Appl. Genet. Mol. \textbf{6}, 1--21 (2007)

\bibitem{vanderlaan2008estimation} van der Laan, M.J: Estimation based on case-control designs with known prevalence probability. Int. J. Biostat. \textbf{4} (2008)

\bibitem{vanderlaan2011targeted} van der Laan, M.J., Rose, S.: Targeted Learning: Causal Inference for Observational and Experimental Data. Springer, New York (2011)

\bibitem{vanderlaan2013estimating} van der Laan, M.J., Petersen, M., Zheng, W.: Estimating the effect of a community-based intervention with two communities. J. Causal Inf. \textbf{1}, 83--106 (2013)

\bibitem{vanderlaan2014causal} van der Laan, M.J.: Causal inference for a population of causally connected units. J. Causal Inf. \textbf{2}, 13--74 (2014)

\bibitem{vanderlaan2014targeted2} van der Laan, M.J.: Targeted estimation of nuisance parameters to obtain valid statistical inference. Int. J. Biostat. \textbf{10}, 29--57 (2014)

\bibitem{vanderlaan2014targeted} van der Laan, M.J.: Targeted learning: From MLE to TMLE. In: Lin, X., Genest, C., Banks, D.L., et al. (eds.) Past, Present, and Future of Statistical Science, pp. 465-480. Chapman \& Hall, London (2014)

\bibitem{vandervaart1996weak} van der Vaart, A.W., Wellner, J.A.: Weak Convergence and Empirical Processes. Springer, New York (1996)

\bibitem{vandervaart2000asymptotic} van der Vaart, A.W.: Asymptotic Statistics. Cambridge University Press (2000)

\bibitem{vandervaart2002semiparametric} van der Vaart, A.W.: Part III: Semiparametric Statistics. In: Bernard, P. (ed.) Lectures on Probability Theory and Statistics, pp. 331-457. Springer, New York (2002)

\bibitem{vandervaart2014higher} van der Vaart, A.W.: Higher order tangent spaces and influence functions. Stat, Sci. \textbf{29}, 679--686 (2014)

\bibitem{vanderweele2009concerning} VanderWeele, T. J.: Concerning the consistency assumption in causal inference. Epid. \textbf{20}, 880--883 (2009)

\bibitem{vanderweele2011weighting} VanderWeele, T.J., Vansteelandt, S.: A weighting approach to causal effects and additive interaction in case-control studies: marginal structural linear odds models. Am. J. Epid. \textbf{174}, 1197--1203 (2011)

\bibitem{vanderweele2014invited} VanderWeele, T.J., Vansteelandt, S.: Invited commentary: some advantages of the relative excess risk due to interaction (RERI) - towards better estimators of additive interaction. Am. J. Epid. \textbf{179}, 670--671 (2014)

\bibitem{vanderweele2015explanation} VanderWeele, T.J.: Explanation in Causal Inference: Methods for Mediation and Interaction. Oxford University Press (2015)

\bibitem{zheng2010asymptotic} Zheng, W., van der Laan, M.J.: Asymptotic theory for cross-validated targeted maximum likelihood estimation. U.C. Berkeley Division of Biostatistics Working Paper Series. \textbf{273}, 1--58 (2010)

\end{thebibliography}
\end{document}